\documentclass[12pt]{article}
\usepackage[final]{epsfig}	
\usepackage{amsfonts,amssymb,latexsym,
amscd,theorem,epsfig,graphics}
\newtheorem{theorem}{Theorem}
\newtheorem{lemma}{Lemma}[section]

\newtheorem{remark}[lemma]{Remark}

\newtheorem{corollary}[theorem]{Corollary}

\begin{document}
\newcommand{\eps}{{\varepsilon}}
\newcommand{\proofend}{$\Box$\bigskip}
\newcommand{\C}{{\mathbf C}}
\newcommand{\Q}{{\mathbf Q}}
\newcommand{\R}{{\mathbf R}}
\newcommand{\Z}{{\mathbf Z}}
\newcommand{\RP}{{\mathbf {RP}}}

\title {Remarks on magnetic flows and magnetic billiards,  Finsler metrics and a magnetic  analog of
Hilbert's fourth problem}
\author{Serge Tabachnikov \\
{\it Department of Mathematics, Penn State University}\\
{\it University Park, PA 16802, USA}\\
 {\it tabachni@math.psu.edu} 
}
\date{}
\maketitle

\begin{abstract} We interpret magnetic billiards as Finsler ones and describe an analog of the string
construction for magnetic billiards. Finsler billiards for which the law ``angle of incidence equals
angle of reflection" are described. We characterize the Finsler metrics in the plane whose geodesics
are circles of a fixed radius. This is a magnetic analog of Hilbert's fourth problem asking to
describe the Finsler metrics whose geodesics are straight lines.

{\it MSC}: 53B40, 53C60, 53C65, 70F35

{\it Key words}: magnetic flows, magnetic billiards, Finsler metrics, Finsler billiards, Hilbert's fourth
problem, Pompeiu problem, integral geometry
\end{abstract}
\bigskip

\section{Introduction and background material}

This paper concerns the motion of a charged particle in a magnetic field, a popular object 
of study in mathematics and mathematical physics. In the Euclidean plane, the strength of the
magnetic field is given by a function $B(x_1,x_2)$, and the particle moves with a constant speed,
satisfying the equation
\begin{eqnarray} 
\ddot x =  B(x_1,x_2) J \dot x\quad {\rm where} \quad J(v_1, v_2) = (-v_2, v_1) \label{flow}
\end{eqnarray}
If one fixes the speed $|v|$ then the magnetic field prescribes the curvature of the trajectory at every
point. In particular, if the field is constant then the trajectories are circles of the Larmor radius
$|v|/|B|$. Our sign convention is that if $B > 0$ then the circles are traversed in the
counterclockwise direction. 

In general, a magnetic field on a Riemannian manifold $M$ is a closed differential 2-form $\beta$,
and the magnetic flow is the Hamiltonian flow of the Riemannian Hamiltonian function
$|p|^2/2$ on the cotangent bundle $T^* M$ with respect to the twisted symplectic structure $\omega +
\pi^*(\beta)$ where $\omega= dp \wedge dq$ is the standard symplectic structure  on $T^* M$ and $\pi:
T^* M \to M$ is the projection. We refer to \cite{Ar1,Bi, Gi, Le, P-S} for a variety of results on
magnetic flows on Riemannian manifolds.

If a charged particle is confined to a domain with ideally reflecting boundary then one has a magnetic
billiard. The particle moves inside according to equation  (\ref{flow}) and undergoes elastic
reflections off the boundary: the tangential component of the velocity remains the same and the normal
one changes sign. In dimension two, this amounts to the familiar law of geometrical optics: the angle of
incidence equals that of reflection. Magnetic billiards have attracted a considerable attention: see
\cite {Ber1,Ber2,B-R,Gu,Tas1,Tas2,Tas3,Zh}; see also \cite{Tab} for a survey of various aspects of 
billiard systems. 

In this paper we interpret a magnetic flow as a geodesic flow of a Finsler metric. We mostly consider
the 2-dimensional case. In Section 2 we interpret magnetic billiards as Finsler ones and describe
 the magnetic version of the string construction that recovers a billiard table by a caustic of the
billiard map. We also characterize the Finsler metrics for which the Finsler billiard enjoys the
familiar law ``angle of incidence equals angle of reflection" (Theorem \ref{converse} and Corollary
\ref{descr}). In Section 3 we describe the Finsler metrics in the plane whose geodesics are circles of a
fixed radius; we give analytic and synthetic descriptions in Theorems \ref{anal},\ref{str},\ref{syn}.
This  is an analog of the celebrated Hilbert's fourth problem of describing the Finsler metrics in a
domain in projective space whose geodesics are straight lines \cite{Al1,Al3,Bu,Po}. Our solution has an
unexpected connection to another classical and well studied question: the Pompeiu problem
\cite{Za1,Za2,Za3}.
\smallskip

We will now review basics of Finsler geometry (see, e.g., \cite{Al2,Ar3,B-C-S,Ru}) and Finsler
billiards, recently introduced in \cite{G-T}.  Finsler geometry describes the propagation of light in an
inhomogeneous  anisotropic medium. This means that the velocity of light depends on the point and the
direction. There are two equivalent descriptions of this process corresponding to the Lagrangian and
the Hamiltonian approaches in classical mechanics, and we will mostly use the former.

The optical properties of a medium are described by a quadratically convex smooth hypersurface, called
the indicatrix, in the tangent space at each point. The indicatrix consists of the velocity vectors of
the  propagation of light at a point in all directions. It plays the role of the unit sphere in
Riemannian geometry. We do not assume that the indicatrices are centrally symmetric.

Equivalently, a Finsler metric on a manifold $M$ is determined by a smooth nonnegative fiber-wise convex
Lagrangian function $L(x,v)$ on the tangent bundle $TM$, homogeneous of degree 1 in the velocity:
$$L(x,tv) = tL(x,v)\quad {\rm for\  all}\quad t>0.$$
The restriction of $L$ to a tangent space $T_xM$ gives the Finsler length of vectors in $T_xM$, and
the indicatrix at $x$ is the unit level hypersurface of $L(x,v)$.  Given a smooth curve $\gamma:[a,b]
\to M$, its length is 
$$ {\cal L} (\gamma) = \int_a^b L(\gamma(t), \gamma'(t))\ dt. $$
The integral does not depend on the parameterization. A Finsler geodesic is an extremal of the
functional ${\cal L}$. The Finsler geodesic flow is the flow in $TM$ in which the foot point of a vector 
moves along the Finsler geodesic tangent to it, so that the vector remains tangent to this geodesic and
preserves its Finsler length. The Finsler geodesic flow is described by the Euler-Lagrange equation
\begin{eqnarray} 
 dL_v(x,v)/dt - L_x(x,v)=0 \quad {\rm or} \quad L_{vv}{\dot v} + L_{vx} v - L_x=0.  \label{EL}
\end{eqnarray} 

The dual, Hamiltonian approach describes the propagation of light in terms of wave fronts and the
Finsler geodesic flow as a Hamiltonian flow in the cotangent bundle $T^*M$. Let $I \subset T_x M$ be the
indicatrix. The figuratrix $J \subset T^*_x M$ is the dual hypersurface constructed as follows. Given a
vector $u \in I$, the respective covector $p \in J$ is defined by the conditions:  
$${\rm Ker}\ p = T_u I \quad {\rm and} \quad  p(u) =1.$$
This gives a diffeomorphism $I \to J$, called
the Legendre transform. In the same way as the field of indicatrices determines the Lagrangian $L$, the
field of figuratrices determines a Hamiltonian $H$ on $T^*M$. The  Hamiltonian vector field of the
function  $H$ is also called the Finsler geodesic flow; the Legendre transform identifies the two flows.
\smallskip

\noindent {\bf Example 1: Hilbert's fourth problem in dimension two.} The Euclidean metric
is given by the Lagrangian $L(x,v) = |v|$; its geodesics are straight lines. Such metrics are called
projective. Following  \cite{Al1}, let us describe all symmetric projective Finsler
metrics in the plane, that is, a solution to Hilbert's fourth problem in dimension 2. 

A synthetic approach, due to Busemann,  makes use of integral geometry, namely, the Crofton formula
\cite{Sa}. Consider the set of oriented lines in the plane, topologically, the cylinder. An oriented
line can be characterized by its direction $\alpha \in [0, 2\pi)$ and its signed distance $p$ from the
origin.  The 2-form $\omega_0 =  dp \wedge d\alpha$ is the standard area form on the space of oriented
lines; this symplectic form is a particular case of a symplectic structure on the space of trajectories
of a Hamiltonian system on a fixed energy level, in particular, the space of oriented geodesics of a
Finsler metric -- see, e.g., \cite {A-G} and a discussion in Section 3. The Crofton formula gives the
Euclidean length of a plane curve $\gamma$ in terms of $\omega_0$. The curve determines a function on
the space of oriented lines, the number of intersections of a line $l$ with $\gamma$. Then
\begin{eqnarray} 
{\rm length}(\gamma) = (1/4) \int \#(l \cap \gamma)\ \omega_0. \label{Croft}
\end{eqnarray}  
Let $f(p,\alpha)$ be a positive continuous function. Then $\omega = f(p,\alpha)\ dp \wedge d\alpha$ is
also an area form on the space of oriented lines. Formula (\ref{Croft}), with $\omega$  replacing
$\omega_0$, defines a projective Finsler metric, and all such metrics can be obtained by an
appropriate choice of the function $f$. 

Next we describe an analytic solution to Hilbert's fourth problem in dimension two. First, a
Lagrangian $L(x,v)$, homogeneous of degree 1 in $v$, gives a projective Finsler metric if and only if
the mixed second partial derivative matrix $L_{xv}$ is symmetric; this is Hamel's theorem of 1903, and it
holds in any dimension.  The Lagrangians satisfying Hamel's condition have the following integral
representation:
\begin{eqnarray} 
L(x_1,x_2,v_1,v_2)=\int_0^{2\pi} |v_1 \cos \phi + v_2 \sin \phi|\ f(x_1 \cos \phi + x_2 \sin
\phi, \phi)\ d\phi \label{Hamel}
\end{eqnarray}
where $f(p,\phi)$ is a smooth positive function on the cylinder representing the space of oriented
lines. Moreover, if $f$ is even in $\phi$  then it is uniquely determined by $L$. The function $f$
is the same as in (\ref{Croft}): the length of a curve with respect to the Finsler metric (\ref{Hamel})
is given, up to a multiplicative constant, by (\ref{Croft}). If $f$ depends on the angle $\phi$ only
then one obtains a translation invariant metric, called a Minkowski metric. If $f$ is a constant then
one has the Euclidean metric.
\medskip

Let $M$ be a 2-dimensional Finsler manifold with boundary, a curve $N$. The Finsler billiard system
is defined in \cite{G-T} as follows. A point moves inside $M$ freely,
according to the  Finsler geodesic flow, until it hits the boundary. The reflection is described in
terms of the indicatrix $I$ at the impact point $x$ -- see figure 1. The vectors $u$ and $v$ are the
Finsler unit vectors along the incoming and outgoing trajectories. The tangent lines to $I$ at $u$ and
$v$ are concurrent with the tangent line to $N$ at $x$.
This definition satisfies a variational principle: for every points $a, b \in
M$, the reflection point $x$ extremizes the Finsler length $|ax| + |xb|$. If the indicatrix is a circle
centered at the origin then the vectors $u$ and $v$ make equal angles with the boundary curve $N$; this
is the familiar law of Euclidean billiard reflection. The multi-dimensional version of the Finsler
billiard reflection is defined similarly, and we do not dwell on it -- see \cite{G-T}.
\smallskip

\begin{figure}[ht]
\centerline{\epsfbox{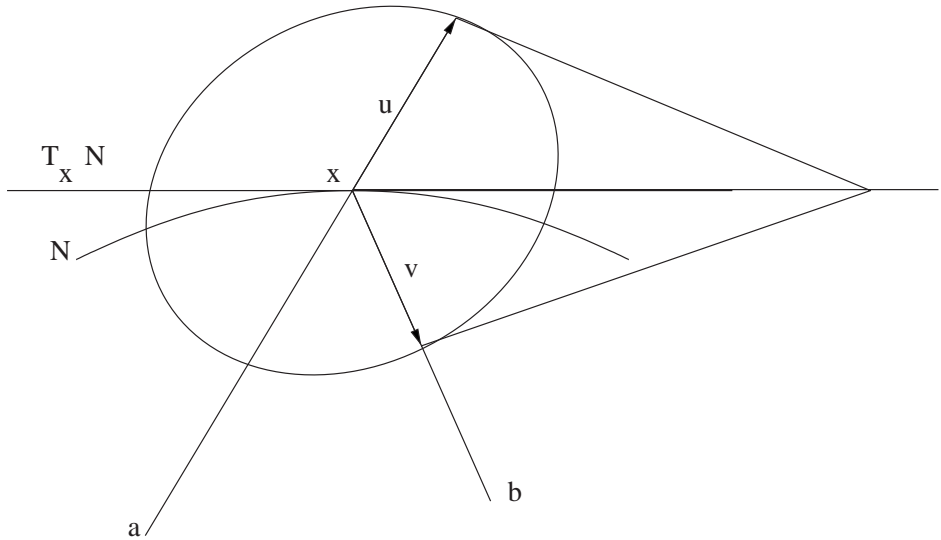}}
\caption{Finsler billiard reflection}
\end{figure}

\noindent {\bf Example 2:  projective Finsler billiard reflection.} Consider a symmetric projective
Finsler metric (\ref{Hamel}). In polar coordinates, $v_1= r \cos \alpha, v_2 = r \sin \alpha$, and
\begin{eqnarray}
L(x_1, x_2, r, \alpha)/2 = r \int_{\alpha - \pi/2}^{\alpha + \pi/2} \cos (\alpha - \phi)\ f(x_1 \cos
\phi + x_2 \sin \phi, \phi)\ d\phi. \label{polar}
\end{eqnarray}
Let $\alpha$ be the direction of the billiard curve at the impact point $x$ and $\beta$ and $\gamma$
the directions of the incoming and the outgoing billiard trajectories -- see figure 2. The projective
Finsler reflection law specializes to the following formula.

\begin{lemma} \label{projlaw}  One has:
$$
\int_{\gamma - \pi/2}^{\beta - \pi/2} \cos (\alpha - \phi)\ f(x_1 \cos \phi
+ x_2 \sin \phi, \phi)\ d\phi =$$
$$ \int_{\gamma + \pi/2}^{\beta + \pi/2} \cos (\alpha - \phi)\ f(x_1
\cos \phi + x_2 \sin \phi, \phi)\ d\phi. 
$$
\end{lemma} 

\begin{figure}[ht]
\centerline{\epsfbox{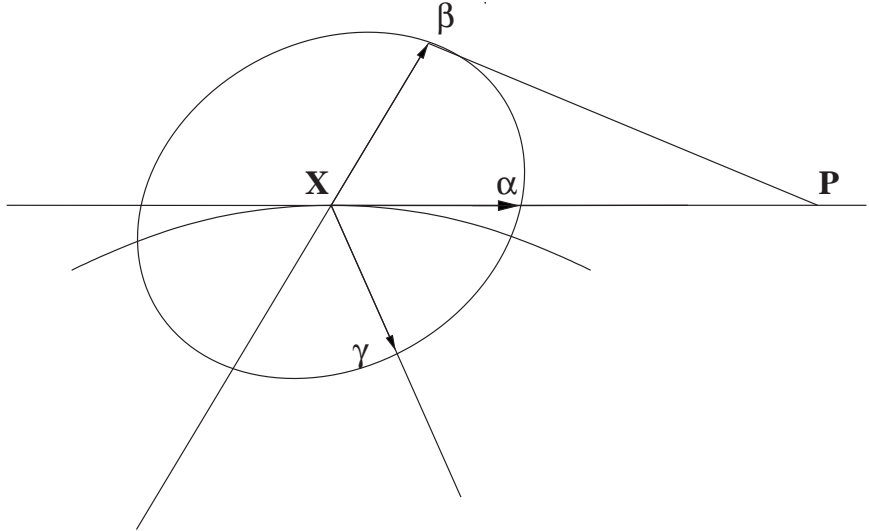}}
\caption{Deriving the projective Finsler billiard reflection law}
\end{figure}

For example, if $f=1$ then integration yields: $\cos (\beta - \alpha) = \cos (\gamma -
\alpha)$ or $\beta - \alpha = \gamma -\alpha$, the familiar law of equal angles. 
\smallskip

{\bf Proof}.
Denote the integral in (\ref{polar}) by $g(x,\alpha)$. Then the polar equation of the indicatrix at
point $x$, chosen as the origin, is $r=1/g(x,\alpha)$. It is a matter of a straightforward calculation
to find the coordinates of the intersection point $P$ in figure 2:
$$P(\beta) = \frac{ (\cos \alpha, \sin \alpha)}{g(\beta) \cos (\alpha - \beta) + g'(\beta) \sin (\alpha -
\beta)}.$$ 
Equating $P(\beta)$ and $P(\gamma)$ yields:
\begin{eqnarray}
g(\beta) \cos (\alpha - \beta) + g'(\beta) \sin (\alpha - \beta) = g(\gamma) \cos (\alpha - \gamma) +
g'(\gamma) \sin (\alpha - \gamma). \label{eqnn}
\end{eqnarray}
It follows from (\ref{polar}) that
$$g'(\beta) = - \int_{\beta - \pi/2}^{\beta + \pi/2} \sin (\beta - \phi)\ f(x_1 \cos \phi
+ x_2 \sin \phi, \phi)\ d\phi,$$
and similarly for $\gamma$. It remains to substitute to (\ref{eqnn}) and to collect terms.
\proofend

\section{Magnetic billiards as Finsler billiards}

Consider the plane motion of a charged particle in a magnetic field with strength $B(x_1, x_2)$. The
Lagrangian for this motion is
$$\bar L(x,v) = \frac{1}{2} |v|^2 + \alpha(x) (v)$$
where $\alpha(x) = f(x_1, x_2)\ dx_1 + g(x_1, x_2)\ dx_2$ is a differential 1-form such that $d \alpha =
-B(x_1, x_2)\ dx_1 \wedge dx_2$. The choice of $\alpha$ is not unique: one can always add a closed
1-form to a Lagrangian without effecting the dynamics. The Euler-Lagrange equation for
$\bar L$ is (\ref{flow}). In particular, the Lagrangian for a constant magnetic field is  
$$\bar L(x,v) = \frac{1}{2} |v|^2 + \frac{B}{2} [v,x]$$
where $[\ ,\ ]$ is the cross-product.

Following the Maupertuis principle (see, e.g., \cite{Ar3}), we replace the Lagrangian $\bar L$ by
\begin{eqnarray} 
L(x,v) = |v| + \alpha(x) (v). \label{Lagr}
\end{eqnarray}
The extremals of the Lagrangian (\ref{Lagr}) coincide with those of $\bar L$, corresponding to the
motion with the unit speed. In particular, the extremals of
\begin{eqnarray}
L(x,v) = |v| + \frac{1}{2R} [v,x] \label{const}
\end{eqnarray}
are the counterclockwise oriented circles of radius $R$.

The Lagrangian (\ref{Lagr}) defines a non-symmetric Finsler metric in the domain where $L(x,v) >0$ for
all $v \neq 0$. This is the case if $|\alpha(x)| < 1$, and we assume this condition to hold throughout
this section. In other words, we assume that the magnetic field is sufficiently weak. Under this
assumption, we consider the unit speed magnetic flow as the Finsler geodesic flow.

Consider a plane domain  and the magnetic billiard inside it. One also has the Finsler billiard inside
the domain, associated with the Lagrangian (\ref{Lagr}). One expects the two systems to coincide, that
is, to have the same reflection laws. 

\begin{theorem} \label{equal} The Finsler billiard reflection law, associated with the Lagrangian
(\ref{Lagr}), is the law of equal angles: the angle of incidence equals the angle of reflection. 
\end{theorem}

{\bf Proof.} The indicatrix of the Finsler metric at point $x$ is given by the equation
$|v| + \alpha(x)(v) =1$. Choose Cartesian coordinates in such a way that $\alpha(x)(v) = tv_1$ for some
$t \in \R$. Then the equation of the indicatrix can be rewritten as
\begin{eqnarray}
(1-t^2)^2 \left(v_1 + \frac{t}{1-t^2}\right)^2 + (1-t^2) v_2^2 =1. \label{conic}
\end{eqnarray}
Recall that the equation of a conic, centered at the origin in the $(v_1,v_2)$-plane, is
$$\frac{(v_1+c)^2}{a^2} + \frac{v_2^2}{b^2}=1$$
where $a^2-c^2=b^2$. Clearly, (\ref{conic}) has this form. Hence the indicatrix is an ellipse, centered
at the origin.

\begin{figure}[ht]
\centerline{\epsfbox{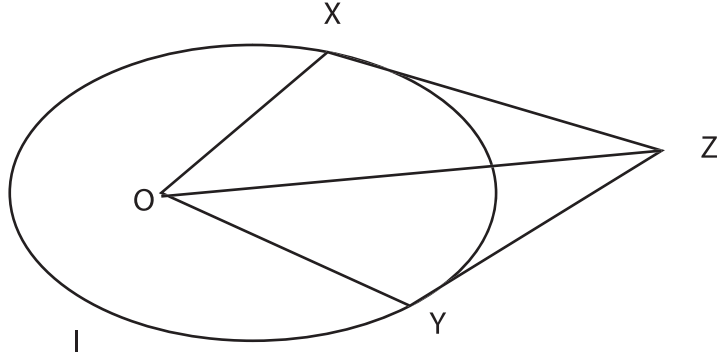}}
\caption{Poncelet's first little theorem}
\end{figure}

Thus the  theorem reduces to the following geometrical property of conics. Let $I$ be a
conic with focus $O$,  let $X, Y \in I$, and let $Z$ be the intersection point of the
tangent lines to $I$ at $X$ and $Y$. Then the line $OZ$ bisects the angle $XOY$ -- see figure 3. This
property  holds indeed; it is known as the Poncelet ``first little theorem", see \cite{Be}. This
completes the proof. \proofend

In fact, the law of equal angles is characteristic of the conics, centered at a focus.

\begin{theorem} \label{converse} Let $I$ be a smooth  plane curve, star-shaped with
respect to  point $O$, with the following property. Let $X, Y \in I$ be arbitrary points,  $Z$ be
the intersection point of the tangent lines to $I$ at $X$ and $Y$; then the line $OZ$ bisects the angle
$XOY$.  It follows that $I$ is a conic with focus $O$.
\end{theorem}

{\bf Proof.} Let $O$ be the origin, and give $I$ a parameterization $I(t)$ so that $[I(t), I'(t)] \equiv
1$. Then $I''(t)=-f(t) I(t)$ for some function $f(t)$; thus we view $I$ as an orbit in a central force
field. We claim that $f(t) = C/|I(t)|^3$ for a constant $C$. Assuming this claim, it follows that $I$
is an orbit in Newton's force field, and therefore a conic with focus $O$.

Let $X=I(t_1)=I_1, Y=I(t_2)=I_2$. A direct computation yields the point $Z$:
$$Z = I_1 + \frac{[I_2 - I_1, I_2']}{[I_1', I_2']} I_1' = I_2 + \frac{[I_1 - I_2, I_1']}{[I_2', I_1']}
I_2'.$$ 
The equal angle condition reads: $[Y,Z]/|Y| = [Z,X]/|X|$, or
\begin{eqnarray}
|I_2|\ (1-[I_1, I_2']) = |I_1|\ (1-[I_2, I_1']). \label{conv}  
\end{eqnarray}
Now set: $t_1=t, t_2=t+\eps$ and shorthand $I(t)$ to $I$ and $f(t)$ to $f$. Then the Taylor expansion
yields:
$$I_2= I \left(1-\frac{\eps^2}{2} f - \frac{\eps^3}{6} f'\right) + I' \left(\eps - \frac{\eps^3}{6}
f\right) +  O(\eps^4),$$
$$I_2'= I'\left(1-\frac{\eps^2}{2} f - \frac{\eps^3}{3} f'\right) - I \left(\eps f +
\frac{\eps^2}{2} f' + \frac{\eps^3}{6}(f''-f^2)\right) + O(\eps^4),$$
and
$$|I_2|=|I| + \eps \frac{I \cdot I'}{|I|} + O(\eps^2).$$
Substitute to (\ref{conv}) and collect terms to obtain:
$|I|^2 f' + 3 I \cdot I' f = 0.$ This differential equation is easily solved: $f'/f= -3 I \cdot
I'/|I|^2$ and hence $f=C/|I|^3$, as claimed. 
\proofend

As a consequence, we obtain a description of Finsler metrics  for which the Finsler
billiard reflection law is the law of equal angles.

\begin{corollary} \label{descr} The Finsler billiard reflection satisfies the law \ ``angle of incidence
equals angle of reflection" for every billiard curve if and only if the metric is given by a Lagrangian
\begin{eqnarray}
L(x,v) = f(x) (|v| + \alpha(x) (v)) \label{conf}
\end{eqnarray}
where $f(x)$ is a non-vanishing function and $\alpha(x)$ is a 1-form.
\end{corollary}

{\bf Proof}. Replacing a metric by a conformally-equivalent one changes the indicatrices by a dilation
and does not effect the law of equal angles. Theorem \ref{equal} implies that the metrics (\ref{conf})
satisfy this law of equal angles. Conversely, if this law  holds then, by Theorem
\ref{converse}, the indicatrices are ellipses (depending on the point of the plane), centered at their
foci. A general equation of such an ellipse is $f(x) (|v| + \alpha(x) (v))=1$, and the result follows.
\proofend

It is shown in \cite{G-T} that some familiar properties of the usual billiards extend to the Finsler
ones. Although \cite{G-T} concerned symmetric Finsler metrics, the results hold in the non-symmetric
case as well; however one should be careful with the order of points: the distance from $A$ to $B$ may
differ from the distance from $B$ to $A$. Let us consider the case of a constant magnetic field, that
is, the Finsler metric given by the Lagrangian (\ref{const}) whose geodesics are counterclockwise oriented
arcs of radius $R$. 

Let $A$ and $B$ be two points on an arc of radius $R$ with the center $C$ and the angle measure
$\theta$. Denote the Finsler distance between points by $d(A,B)$ and identify points with their
position vectors. Let ${\cal L} (\gamma)$ denote the Finsler length of a curve $\gamma$.

\begin{lemma} \label{explicit} One has:
$$d(A,B) = \frac{1}{2} \theta R + \frac{1}{2R} [B-A,C].$$
For a simple oriented closed curve $\gamma$, one has:
\begin{eqnarray}
{\cal L} (\gamma) = l(\gamma) - \frac{1}{R} S(\gamma) \label{length}
\end{eqnarray}
where $l(\gamma)$ and $S(\gamma)$ are the Euclidean length and the Euclidean signed area bounded by
$\gamma$.
\end{lemma}

{\bf Proof}. To obtain (\ref{length}), one integrates $|v| +  [v,x]/(2R)$ over $\gamma$ and makes
use of the fact that $[v,x]/2$ is negative the derivative of the signed area swept by the position
vector of
$\gamma$.

Let $O$ be the origin. The distance $d(A,B)$ equals the integral of $|v| +  [v,x]/(2R)$
over the arc $AB$. The integral of $|v|$ is the arclength of the arc, that is, $\theta R$. The integral
of $[v,x]/(2R)$ equals $-S/R$ where $S$ is the area of the curvilinear triangle $OAB$.  

Assume first that $O=C$. Then the latter area is $\theta R^2/2$, and $d(A,B)=\theta R/2$. If the origin
is translated through vector $C$  then the Lagrangian changes by the term $[v,C]/(2R)$, and its integral
by $[B-A,C]/(2R)$. This yields the first formula.
\proofend

\begin{remark} \label{gen} {\rm Formula (\ref{length}), along with its proof, holds for closed immersed
curves as well: the area term should be understood as the integral of the 1-form $(xdy-ydx)/2$ over the
curve.}
\end{remark}

The orientation of $\gamma$ determines a coorientation: the
pair (coorientation vector, orientation vector) gives the positive orientation of the plane. If $\gamma$
is a counterclockwise oriented simple curve then the positive coorientation is the outward one. Given a
real number $t$, consider the parallel curve $\Gamma(t)$ at distance $t$ from
$\gamma$. The curve $\Gamma(t)$ is the time-$t$ wave front, starting at $\gamma = \Gamma(0)$. More
precisely, one translates the contact elements of $\gamma$ in the orthogonal direction through distance
$t$ (along the coorienting vector, if $t>0$, and in the opposite direction, if $t <0$), and the obtained
1-parameter family of contact elements consists of the contact elements of $\Gamma(t)$. The curve
$\Gamma(t)$ may have singularities, generically, semi-cubic cusps. If $\gamma$ is a
positively oriented circle of radius $r$ then $\Gamma(t)$ is a circle of radius $r+t$.

Formula (\ref{length}) admits the following interpretation.

\begin{lemma} \label{interp} For an oriented simple closed curve $\gamma$, one has:
\begin{eqnarray}
{\cal L} (\gamma) = \frac{1}{R} \left(\pi R^2 - S(\Gamma(-R)) \right) \label{interplength}
\end{eqnarray}
where the area  $S(\Gamma(-R))$ is understood as in Remark \ref{gen}.
\end{lemma}

{\bf Proof}. For a closed immersed curve, the following well known formula holds:
$$S(\Gamma(t)) = S(\gamma) + t l(\gamma) + \pi t^2 w$$
where $w$ is the Whitney winding number of $\gamma$. Therefore the right hand side of (\ref{length})
equals $(1/R)(\pi R^2 - S(\Gamma(-R)))$, and the result follows from Lemma \ref{explicit}.
 \proofend

\begin{remark} \label{dual} {\rm Formula (\ref{interplength}) expresses  Finsler lengths in terms of
areas, and it serves a magnetic analog of the Crofton formula (\ref{Croft}). The curve $\Gamma(-R)$ is
the locus of the centers of positively oriented circles of radius $R$, tangent to $\gamma$ and having
the same orientation as $\gamma$ at the tangency point. The curve $\gamma$ can be reconstructed from
$\Gamma(-R)$ as the envelope of the family of circles of radius $R$ centered at points of $\Gamma(-R)$.  

Formula (\ref{interplength})  also resembles the area-length duality for spherical curves, discussed in
\cite{Ar2,Tab1, Tab2}. Let $\gamma$ be a simple smooth closed curve on the unit sphere, and let $\Gamma$
be its spherically dual curve, namely, the curve $\Gamma(\pi/2)$, in the sense of spherical geometry.
Then one has: $l(\gamma) = 2\pi - S(\Gamma).$}
\end{remark}

Let us return to billiards. Recall that a caustic of a 2-dimensional billiard is a curve $\Gamma$ inside
it with the following property: if a segment of a billiard trajectory is tangent to $\Gamma$ then so is
the reflected segment. Given a convex caustic $\Gamma$, can one reconstruct the billiard table? For the
usual, Euclidean billiard the answer is given by the string construction: a billiard curve $N$ is the
locus of points $X$ of a string of fixed length, wrapped around $\Gamma$ --  \cite{Be, Tab}. 

It is shown in \cite{G-T} that the string construction extends to Finsler billiards as
well. For non-symmetric Finsler metrics one needs to consider oriented caustics $\Gamma$ so that the
orientation of billiard trajectories, tangent to $\Gamma$, agrees with the orientation of $\Gamma$.
Applying these considerations to billiards in a constant magnetic field, we obtain the following
corollary.

\begin{figure}[ht]
\centerline{\epsfbox{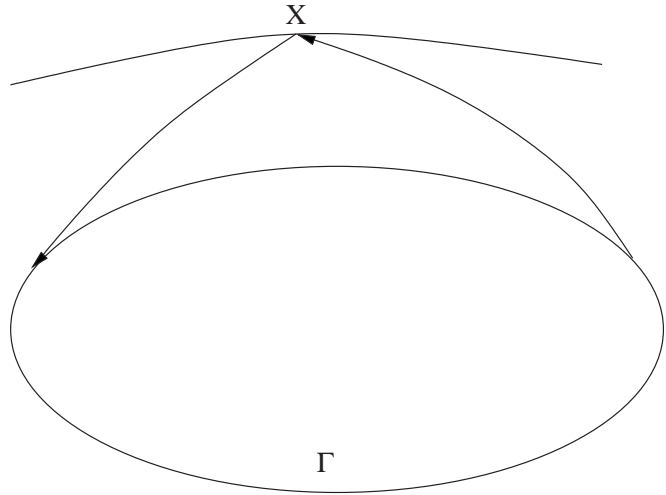}}
\caption{Magnetic string construction}
\end{figure}

Let $\Gamma$ be an oriented closed convex curve. For a point $X$ outside of $\Gamma$, let $F(X)$ be the
Finsler length, associated with the metric (\ref{const}), of the shortest closed  curve from $X$ to $X$
around the obstacle $\Gamma$, whose orientation agrees with that of $\Gamma$ -- see figure 4.

\begin{corollary} \label{string} The level curves of the function $F(X)$ are the boundaries of magnetic
 billiard tables that have  $\Gamma$ as a caustic.
\end{corollary}

Recall the optical property of an ellipse: a ray emanating from one focus reflects to another focus. 
As a particular case of Corollary \ref{string}, one may construct a magnetic analog of an ellipse. 

\begin{corollary} \label{ellipse} Let $A$ and $B$ be fixed points (``foci") and $N$ be the locus of
points $X$ such that $d(A,X) + d(X,B) = {\rm const}$. Then every trajectory of the magnetic billiard,
starting at $A$, reflects in $N$ to $B$.
\end{corollary}

Note that the two foci play different roles and cannot be interchanged in the above formulation. If the
points $A$ and $B$ merge then the ``ellipse" $N$ becomes a Euclidean circle centered at this point.

\section{Finsler metrics whose geodesics are circles of a fixed radius}

In this section we develop a magnetic analog of the solution to Hilbert's fourth problem, outlined in
Example 1 above. Start with an analytic description of the Lagrangians, homogeneous of degree 1
in the velocity, whose extremals are positively oriented circles of radius $R$.

\begin{lemma} \label{desc} The extremals of $L(x,v)$ are positively oriented circles of radius $R$ if
and only if $L$ satisfies the equation:
\begin{eqnarray} 
\frac{|v|}{R} L_{vv} (Jv) + L_{vx} (v) = L_x   \label{eqn}
\end{eqnarray}
where $J(v_1, v_2) = (-v_2, v_1)$. 
\end{lemma}

{\bf Proof}. Let $x(t)$ be a parameterized curve, $v = x'$. The curve is a counterclockwise
oriented circle if and only if
$$\left(\frac{v}{|v|}\right)' = \frac{1}{R} J(v).$$
Differentiate to express the acceleration vector:
\begin{eqnarray}
v' = \frac{|v|}{R} J(v) + \frac{(v \cdot v') v}{|v|^2}. \label{v'}
\end{eqnarray}
Since the Lagrangian is homogeneous of degree 1, the Euler equation 
$L_v v = L$ holds, and hence $L_{vv} (v) = 0.$
It remains to substitute $v'$ from (\ref{v'}) to the Euler-Lagrange equation
$$L_{vv}(v') + L_{vx}(v) = L_x,$$
and the result follows. \proofend

We are ready to prove the main analytical result of this section.

\begin{theorem} \label{anal} Every Lagrangian, homogeneous of degree 1 in the velocity, whose extremals
are positively oriented circles of radius $R$ can be represented, in polar coordinates, as follows:
$$ L(x,v) = L(x_1, x_2, r, \alpha) =  r \Bigl(\int_0^{\alpha+\pi/2} \cos (\alpha - \phi)\ g(x_1+R \cos
\phi, x_2 + R \sin \phi)\ d\phi\ + $$
\begin{eqnarray}
a(x_1, x_2) \cos \alpha + b(x_1, x_2) \sin \alpha\Bigr) \label{main}
\end{eqnarray}
where $g$ is a positive density function in the plane such that the center of mass of every circle of
radius $R$ is its center, and $a, b$ are two functions, satisfying
\begin{eqnarray}
a_{x_2} (x_1, x_2) - b_{x_1} (x_1, x_2) = \frac{1}{R}\ g(x_1 + R, x_2). \label{ini1}
\end{eqnarray}
\end{theorem}

{\bf Proof}. In polar coordinates,
$v_1 = r \cos \alpha,\ v_2 = r \sin \alpha$, and one has: $L(x,v) = |v|\ p(x_1, x_2, \alpha)$ for some
function $p$. 

Fix a point $x=(x_1,x_2)$ and consider the indicatrix $I$ at $x$, chosen as the origin. The polar
equation of $I$ is $r=1/p(x,\alpha)$. Therefore $p(x,\alpha)$ is the support function of the dual curve,
the figuratrix $J$ (see \cite{Sa}). Parameterize $J$ by the angle $\phi$ made by its tangent vector with
the horizontal axis. Let $f(x,\phi)$ be the radius of curvature at point $J(\phi)$ and let $(a(x),
b(x))$ be the coordinates of the point $J(0)$. One has:
$$J'(\phi) = f(x,\phi) (\cos \phi, \sin \phi),$$
and hence
$$J(\alpha + \pi/2) = J(0) + \int_0^{\alpha+\pi/2} f(x,\phi) (\cos \phi, \sin \phi)\ d\phi,$$
see figure 5. It follows that 
$$ p(x,\alpha) = (\cos \alpha, \sin \alpha) \cdot J(\alpha + \pi/2) = $$
\begin{eqnarray}
a(x) \cos \alpha + b(x) \sin \alpha + \int_0^{\alpha+\pi/2} \cos (\alpha - \phi) f(x,\phi)\ d\phi.
\label{rep}
\end{eqnarray}
Differentiating (\ref{rep}) twice, one  recovers the function $f$ from $p$:
$$f(x,\alpha + \pi/2) = p(x,\alpha) + p''(x,\alpha).$$
Note that every function of the form $p + p''$ is $L^2$ orthogonal to $\cos \alpha$ and $\sin \alpha$.
Thus
\begin{eqnarray}
\int_0^{2\pi} f(x,\alpha) \cos \alpha\ d\alpha = \int_0^{2\pi} f(x,\alpha) \sin \alpha\ d\alpha =0;
\label{ort} 
\end{eqnarray}
this also follows from the integral representation (\ref{rep}) and periodicity of $p$ as a
function of $\alpha$.  

\begin{figure}[ht]
\centerline{\epsfbox{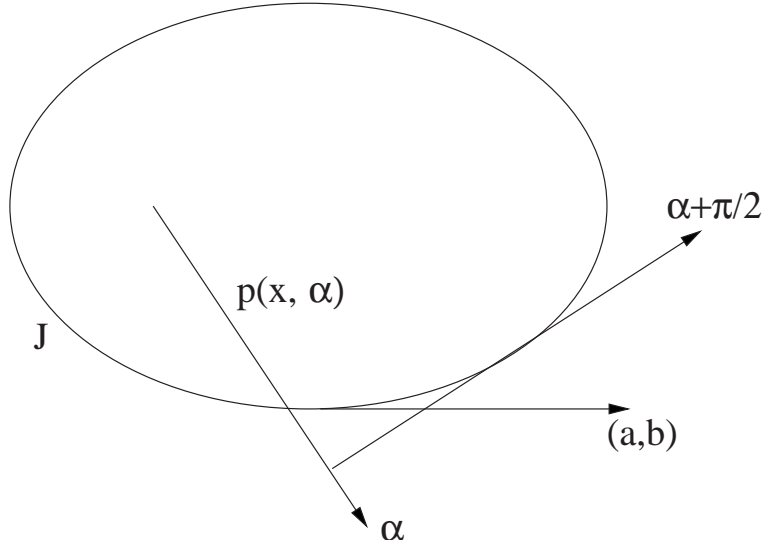}}
\caption{Integral representation of the support function}
\end{figure}

Next we use the equation (\ref{eqn}) in the integral representation (\ref{rep}). One rewrites the
differential operators $\partial^2 v, \partial v \partial x$ and $\partial x$ in polar coordinates and
applies to $L(x_1, x_2, r, \alpha) =  r p(x_1, x_2, \alpha)$, given by (\ref{rep}). Taking into
account that $v = r (\cos \alpha, \sin \alpha)$, a computation reveals that
$$\frac{|v|}{R} L_{vv} (Jv) = \frac{rf(x,\alpha+\pi/2)}{R}\ (-\sin \alpha, \cos \alpha),$$
and
$$L_x - L_{vx} (v) = r \Bigl( \int_0^{\alpha+\pi/2} (\cos \phi\ f_{x_2} (x,\phi) - \sin \phi\
f_{x_1} (x,\phi))\ d\phi\ +$$
$$a_{x_2}(x) - b_{x_1}(x)\Bigr)\ (-\sin \alpha, \cos \alpha).$$
Therefore, by (\ref{eqn}),
$$ 
\frac{1}{R} f(x,\alpha+\pi/2) =  \int_0^{\alpha+\pi/2} (\cos \phi\ f_{x_2} (x,\phi) - \sin \phi\
f_{x_1} (x,\phi))\ d\phi\ + $$
\begin{eqnarray}
a_{x_2}(x) - b_{x_1}(x). \label{int}
\end{eqnarray}
In particular, 
\begin{eqnarray}
\frac{1}{R}\ f(x,0) = a_{x_2}(x) - b_{x_1}(x). \label{ini}
\end{eqnarray}
Differentiating (\ref{int}) with respect to $\alpha$, one gets:
\begin{eqnarray}
\frac{1}{R}\ f_{\alpha} (x,\alpha) = \cos \alpha\ f_{x_2} (x,\alpha) - \sin \alpha\
f_{x_1} (x,\alpha). \label{eqn1}
\end{eqnarray}
We claim that
\begin{eqnarray} 
f(x_1,x_2,\alpha)=g(x_1+R\cos \alpha,x_2+R\sin \alpha)  \label{eqn2}
\end{eqnarray}
for an appropriate function of two variables $g$.

Indeed, consider the vector field 
$$\eta= \frac{1}{R}\ \partial \alpha +\sin \alpha\ \partial x_1 -\cos \alpha\ \partial x_2$$
 on the solid torus
$\R^2 \times S^1$. Then (\ref{eqn1}) can be written as $\eta(f)=0$. The trajectories of $\eta$ are:
$$\alpha(t) = \frac{t}{R},\ x_1(t)= -R\cos\left(\frac{t}{R}\right) +a,\ x_2(t)=
-R\sin\left(\frac{t}{R}\right) +b$$
where $t$ is the ``time" parameter and $a,b$ are constants.
One can take the plane  $\alpha=0$ as a section. Then the $\eta$-invariant function $f$ is
determined by its values on this section, a function of two variables $g$.  Consider a point $(\alpha,
x_1, x_2)$. The trajectory through this point intersects the section at point $(a,b)=(x_1+R\cos \alpha,
x_2+R\sin \alpha)$.  Hence $f(x_1,x_2,\alpha)=g(a,b)$, and (\ref{eqn2}) follows. 

Equations (\ref{ort}) imply that the center of mass of the circle of radius $R$,  with the density
function $g$, centered at $(x_1,x_2)$, is the point $(x_1,x_2)$. Combining (\ref{ini}) and (\ref{eqn2}),
we obtain (\ref{ini1}).
\proofend

\begin{remark} \label{form} {\rm The term $r(a(x) \cos \alpha + b(x) \sin \alpha)$ in (\ref{main}), the
formulation of Theorem \ref{anal}, can be written as $\nu(x) (v)$ where $\nu(x) = a(x)\ dx_1 + b(x)\
dx_2$ is a 1-form. The choice of this form is not unique but $d\nu$ is uniquely determined by the
function
$g$ via (\ref{ini1}). This is consistent with the remark we already made: adding a closed 1-form to the
Lagrangian does not effect the Euler-Lagrange equations.}
\end{remark}

To proceed, we recall basic facts about the symplectic reduction. 
Let $(M,\Omega)$ be a symplectic manifold and $H:M \to \R$ a Hamiltonian function.
Consider the Hamiltonian vector field $\xi = {\rm sgrad}\ H$. Since $H$ is $\xi$-invariant, the field
$\xi$ is tangent to the level hypersurfaces of $H$. Consider such a hypersurface $S$, and assume that
the space of trajectories of $\xi$ on $S$ is a smooth manifold $N = S/\xi$; locally, this is always
the case. The restriction of $\Omega$ to $S$ has a 1-dimensional kernel spanned by $\xi$, and hence
$\Omega|_{S}$ descends to a symplectic structure $\omega$ on $N$. This is the symplectic reduction of
$\Omega$.

One applies this construction as follows. Given a Finsler manifold $M$, the symplectic manifold in
question is the cotangent bundle $T^*M$ with its standard symplectic structure $dp \wedge dx$ where $x
\in M$ is the position and $p \in T^*_x M$ the momentum. The  function $H$ is the Finsler
metric Hamiltonian, and the hypersurface $S$ consists of the unit covectors; it is fibered over $M$ and
the fibers are the figuratrices. The vector field $\xi$ is the Finsler geodesic flow, and the space of
trajectories identifies with the space of non-parameterized oriented geodesics. 

Consider the tangent bundle $TM$ and the unit vector hypersurface $U$ in it. The Legendre transform
$(x,v) \mapsto (x,p=L_v)$ identifies $U$ with $S$ and the Finsler geodesic flow $\zeta$ on $U$ with the
geodesic flow $\xi$ on $S$. The pull-back of the Liouville form $pdx$ is the 1-form $\lambda = L_v dx$ on
$TM$. The form $\lambda$ is a contact form on $U$, and $\zeta$ is its Reeb vector field: $\lambda
(\zeta) = 1,\ i_{\zeta} d\lambda =0$. The reduction of the 2-form $d\lambda$ yields the symplectic
structure on the quotient space $U/\zeta$, the space of oriented Finsler geodesics.

Given a smooth curve $\gamma$ on $M$, one lifts it to the curve $\tilde \gamma$ on $U$ by assigning the
unit tangent vector to every point of $\gamma$. Then the Finsler length of $\gamma$ equals
$$\int_{\tilde \gamma} \lambda.$$
For a reference on this symplectic approach, see, e.g., \cite{Ar3, A-G}.

Now we are in a position to compute the symplectic structure on the space of circles of radius $R$,
associated with the Lagrangian (\ref{main}). A circle is characterized by its center, and the space of
circles is the plane with Cartesian coordinates $(u,v)$.

\begin{theorem} \label{str} The symplectic structure $\omega$ on the space of circles of radius $R$,
associated with the Lagrangian (\ref{main}), is given by the formula:
\begin{eqnarray} 
\omega = - \frac{1}{R}\ g(u,v)\ du \wedge dv. \label{stre}
\end{eqnarray}
\end{theorem}

{\bf Proof}. The manifold $U$ consists of the Finsler unit tangent vectors in the plane
and  has coordinates $\alpha,x_1, x_2$. We use the notation from the proof of Theorem \ref{anal}. The
formulas derived in that proof yield:
$$\zeta = \frac{1}{p(x,\alpha)}\  \left(\cos \alpha\ \partial x_1 + \sin \alpha\ \partial x_2 +
\frac{1}{R}\ \partial \alpha \right),$$
$$\lambda = L_v dx = \left(\int_0^{\alpha + \pi/2} \cos \phi\ f(x,\phi)\ d\phi + a(x) \right)\ dx_1 +$$
$$\left(\int_0^{\alpha + \pi/2} \sin \phi\ f(x,\phi)\ d\phi + b(x) \right)\ dx_2$$
and, taking (\ref{int}) and (\ref{ini}) into account,  
$$d\lambda = f(x, \alpha+\pi/2)\ \left(\cos \alpha\ d\alpha \wedge dx_2 - \sin \alpha\ d\alpha \wedge
dx_1 - \frac{1}{R} dx_1 \wedge dx_2 \right).$$
In view of (\ref{eqn2}),
$$d\lambda = g(x_1-R  \sin \alpha, x_2 + R \cos \alpha)\ \Bigl(\cos \alpha\ d\alpha \wedge dx_2 - $$
\begin{eqnarray}
\sin \alpha\ d\alpha \wedge dx_1 - 
\frac{1}{R} dx_1 \wedge dx_2 \Bigr). \label{forma}
\end{eqnarray}
Now consider the projection $U \to U/\zeta = \R^2$. To compute the symplectic structure $\omega$ in
$\R^2$, consider a section $j: \R^2 \to U$ and let $\omega$ be the pull-back of $d\lambda$; the result
is independent of the choice of $j$. As a section one may take
$$j(u,v) = (\alpha, x_1, x_2)\ \ {\rm with}\ \ \alpha =\frac{\pi}{2},\ x_1=u+R,\ x_2=v,$$
see figure 6. It remains to substitute to (\ref{forma}), and the result follows. 
\proofend

\begin{figure}[ht]
\centerline{\epsfbox{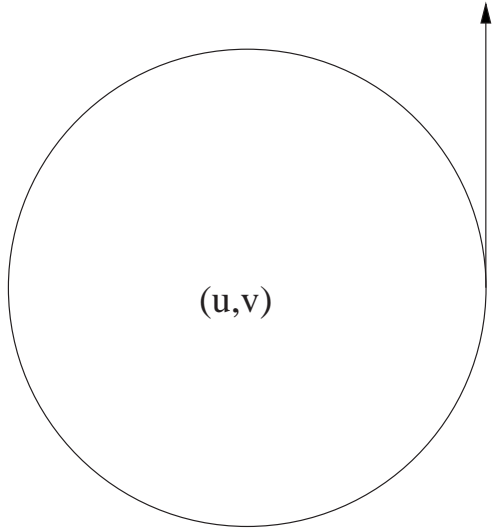}}
\caption{Computing the area form on the space of circles}
\end{figure}

As a consequence, the Finsler metric (\ref{main}) can be recovered, up to summation with a closed 1-form,
from the area form (\ref{stre}) -- see Remark \ref{form}.

Next we consider an analog of formula (\ref{interplength}) for a general Finsler metric
(\ref{main}) whose geodesics are circles of radius $R$. 
The following result expresses the Finsler length in terms of the area form on the space of circles and
is analogous to the synthetic solution to Hilbert's fourth problem, that is, the Crofton formula, used
as a definition of a projective metric.

\begin{theorem} \label{syn} Given an oriented simple closed curve $\gamma$, one has:
\begin{eqnarray}
{\cal L} (\gamma) = S(\Gamma(-R)) + C \label{fla}
\end{eqnarray}
where $S(\Gamma(-R))$ is the area bounded by the curve $\Gamma(-R)$ with respect to the area form
(\ref{stre}), and $C$ is the common Finsler length of all positively oriented circles of radius $R$.
\end{theorem}

{\bf Proof}. Note first that the geodesics are extremals of the length functional ${\cal L}$. The space
of geodesics identifies with the plane which is a critical manifold of ${\cal L}$. A function is constant
on its critical manifold, hence all positively oriented circles of radius $R$ have equal Finsler length.

To prove the result we consider a variation of the curve $\gamma$ and show that  both
sides of (\ref{fla}) have the same variations. This  being established, one can deform $\gamma$ to a
circle of radius $R$ for which the result holds.

Assume that $\gamma$ is parameterized by the Euclidean 
arc-length. Then $\gamma' = (\cos \alpha, \sin \alpha)$. The lift of $\gamma$ to $U$ is the curve
$\tilde \gamma = (\gamma, \alpha)$. Consider a variation of the curve, that is,  a vector field $w$ along
$\gamma$. It is straightforward to compute that the respective variation of $\tilde \gamma$ is the vector
field 
\begin{eqnarray}
\tilde w = w + [\gamma', v']\ \partial \alpha. \label{til}
\end{eqnarray}
 One has:
$${\cal L} (\gamma) = \int_{\tilde \gamma} \lambda$$ 
where $\lambda$ is the contact form as in the proof of Theorem \ref{str}. Therefore the variation of the
length ${\cal L} (\gamma)$ is given by the formula
$$
\int_{\tilde \gamma} i_{\tilde w} d\lambda 
$$
where $d\lambda$ is as in (\ref{forma}). Let $k(t)$ be the curvature at $\gamma(t)$. Then $d\alpha = k
dt$.  A  computation using (\ref{forma}) and (\ref{til}) reveals that
\begin{eqnarray}
\int_{\tilde \gamma} i_{\tilde w} d\lambda = \int g(x_1-R  \sin \alpha, x_2 + R \cos \alpha)\ [\gamma',
w]\ \left(\frac{1}{R} - k \right) dt. \label{var}
\end{eqnarray}

On the other hand, one has: 
$$\Gamma(-R) = (X_1,X_2) = (x_1 - R \sin \alpha, x_2 + R \cos \alpha).$$
Therefore the variation of $\Gamma(-R)$ is given by the vector field 
$$u = (w_1 - R [\gamma', v'] \cos \alpha,\ w_2 - R [\gamma', v'] \sin \alpha).$$
Then
$$dX_1 = (1-Rk) \cos \alpha\ dt,\ dX_2 = (1-Rk) \sin \alpha\ dt.$$
Since $\omega = (-1/R)\ g(X_1, X_2)\ dX_1 \wedge dX_2$, it is straightforward to compute
the variation of the area $S(\Gamma(-R))$:
$$
 \int_{\Gamma(-R)} i_{u} \omega = \int g(x_1-R  \sin \alpha, x_2 + R \cos \alpha)\ [\gamma',
w]\ \left(\frac{1}{R} - k \right) dt. 
$$
This is the same as (\ref{var}), and we are done. \proofend

Note the following corollary of formula (\ref{fla}).

\begin{corollary} \label{equall} The integral of the area form (\ref{stre}) is the same over all discs of
radius $R$.
\end{corollary}

{\bf Proof}. Let $\gamma$ degenerate to a point in (\ref{fla}), so that ${\cal L} (\gamma) =0$. Then
$\Gamma(-R)$ is a circle of radius $R$, and the $\omega$-area, bounded by it, equals $-C$.
\proofend

\begin{remark} \label{alt} {\rm One can give a somewhat different proof of Theorem \ref{syn}
that does not use the specifics of the Euclidean plane and applies to other surfaces, for example, the
sphere. Let us outline the argument. Pick a point $O$ inside $\gamma$ and consider an infinitesimally
small loop
$\delta$ around $O$ whose orientation is the same as that of $\gamma$. The Finsler unit tangent
vector fields to $\gamma$ and to $\delta$ extend to a unit vector field in the annulus $A$ bounded by
$\gamma$ and $\delta$. This vector field provides a lift $\tilde A$ of the annulus to $U$, and
$\partial \tilde A = \tilde \gamma - \tilde \delta$. By Stokes' theorem,
$$\int_{\tilde \gamma} \lambda - \int_{\tilde \delta} \lambda = \int_{\tilde A} d\lambda.$$
The second integral on the left hand side is infinitesimally small. The integral on the
right can be understood as the symplectic area of the set of circles of radius $R$ whose centers lie
between the curves $\Gamma(-R)$ and the circle of radius $R$, centered at $O$, and (\ref{fla})
follows. This also shows that the symplectic area of a circle of radius $R$ is independent on its
choice. }
\end{remark}

One can revert the arguments and and define the respective Finsler metric, as in
Theorem \ref{anal}, starting with an area form $\omega = g(x_1,x_2)\ dx_1 \wedge dx_2$,
satisfying the property that the $\omega$-area of every disc of radius $R$ is the same. Then the
function $g(x_1,x_2)$ should be $L^2$ orthogonal to cosine and sine on every circle of radius $R$.

\begin{lemma} \label{equi} The integrals of a function $g$ over all
discs of radius $R$ is the same if and only if $g$ is orthogonal to cosine and sine on every circle of
radius $R$.
\end{lemma}

{\bf Proof}. Let $S$ be a circle of radius $R$ with center $x=(x_1,x_2)$. Consider its  variation  given
by an infinitesimal parallel translation through vector $v=(v_1,v_2)$. The variation of the
$\omega$-area of the disc is
$$\int_S i_v \omega = \int_0^{2\pi} g(x_1 + R\cos \alpha, x_2 + R \sin \alpha)\ (v_1 \sin \alpha - v_2
\cos \alpha)\ d\alpha.$$
This vanishes for all $v$ if and only if $g$ is orthogonal to $\cos \alpha$ and  $\sin \alpha$. 
\proofend

How restrictive are these two equivalent conditions on function $g$? This question goes to the heart of
the Pompeiu problem, see \cite{Za1,Za2,Za3}. Given a compact set $K$, one considers the
continuous functions with zero integrals over all isometric images of $K$. For which sets $K$ must such
functions  be identically zero? D. Pompeiu, who posed this problem in the late 1920-s, erroneously
thought that the disc in the plane has this property. In fact, if $K$ is a disc then there are plenty of 
functions with zero integrals over all congruent discs; although there is a wealth of results on this
subject, the general solution to the Pompeiu problem is not known yet.

If a function $g$ has a constant integral over all discs of radius $R$ then it can be written as Const +
$h(x_1,x_2)$ where $h$ has zero integrals over all discs of radius $R$.
The following result, standard in the literature on the Pompeiu problem, provides a substantial supply of
such functions. We need to recall the definition of the Bessel functions.

The Bessel functions $J_n(w),\ n \in \Z$, are defined by the generating function
\begin{eqnarray}
\exp \left(\frac{w}{2} \left( t - \frac{1}{t} \right) \right) = \sum_{n=-\infty}^{\infty} J_n(w)\ t^n. 
\label{Bes}
\end{eqnarray}
An explicit formula is as follows:
$$J_n(w) = \sum_{j=0}^{\infty} \frac{(-1)^j w^{n+2j}}{2^{n+2j}(n+j)!j!},\quad n\geq 0;\quad J_{-n}(w) =
J_n(-w).$$
 We will need the following property: 
\begin{eqnarray}
\int J_0 (w)\ w dw = w J_1 (w).\label{Bprop}
\end{eqnarray}

\begin{lemma} \label{pomp} Let $a$ be a root of the first Bessel function $J_1$, and let $f(\beta)$ be
a function on the circle. Then the functions
\begin{eqnarray}
h(x_1,x_2) = \int_0^{2\pi} \cos \left(\frac{a}{R} \left(x_1 \cos \beta + x_2 \sin \beta\right)\right)\
f(\beta)\ d\beta 
\label{pompform}
\end{eqnarray}
and
$$ h(x_1,x_2) = \int_0^{2\pi} \sin \left(\frac{a}{R} \left(x_1 \cos \beta + x_2 \sin
\beta\right)\right)\ f(\beta)\ d\beta $$
have zero integrals over all discs of radius $R$.
\end{lemma}

One may also take linear combinations of such function over different roots of $J_1$.
\smallskip

{\bf Proof}. Let $D$ be the disc of radius $R$ centered at the origin, and let $\xi$ be its
characteristic function.  The condition on function $h$ reads: $\xi * h = 0$ where $*$ denote the
convolution. Take the Fourier transform to obtain: 
\begin{eqnarray}
\hat \xi\  \hat h = 0. \label{FT}
\end{eqnarray}
Let us compute $\hat \xi$:
$$\hat \xi (\lambda) = \int_D e^{-i\lambda x}\ dx = \int_0^R \int_0^{2\pi} e^{-ir\rho \cos(\alpha -
\beta)}\ d\alpha\ r dr$$
where $x=r(\cos \alpha, \sin \alpha), \lambda = \rho(\cos \beta,\sin \beta)$. One has:
$$-ir\rho \cos(\alpha - \beta) = \frac{r\rho}{2} \left(e^{i\theta} - e^{-i\theta}\right)$$
where $\theta = \alpha - \beta - \pi/2$. Using the definition of Bessel functions (\ref{Bes}), it
follows that
$$\int_0^{2\pi} e^{-ir\rho \cos(\alpha - \beta)}\ d\alpha = \sum_n J_n(r\rho) \int_0^{2\pi}
e^{in\theta}\ d\theta = 2\pi J_0 (r\rho).$$
By (\ref{Bprop}), one has:
$$\int_0^R J_0 (r\rho) rdr  = \frac{R}{\rho} J_1 (R\rho)$$
and hence
$$\hat \xi (\lambda) =  \frac{2\pi R}{\rho} J_1 (R|\lambda|).$$
The condition (\ref{FT})  holds if the support of $\hat h$ is contained in the union of circles,
centered at the origin, whose radii are $a/R$ where $a$ is a root of $J_1$. Fix one such root and let
$$\hat h (\rho, \beta) = f(\beta)\ \delta_{\frac{a}{R}} (\rho).$$ 
Taking the inverse Fourier transform yields:
$$h(x_1,x_2) = \int_0^{2\pi} \int_0^{\infty} f(\beta) \delta_{\frac{a}{R}} (\rho) e^{i \rho (x_1 \cos
\beta + x_2 \sin \beta)}\ \rho d\rho\ d\beta = $$
$$\frac{a}{R} \int_0^{2\pi} f(\beta) e^{i \frac{a}{R} (x_1 \cos \beta + x_2 \sin \beta)}\ d\beta.$$
One concludes by taking the real and imaginary parts.
\proofend

For example, let $f$ in (\ref{pompform}) be the delta function $\delta_0$. Then $h= \cos(ax_1/R)$.
Substitute $g=1+ \cos(ax_1/R)$ into (\ref{main}) to obtain an ``exotic" Finsler metric whose geodesics
are circles of radius $R$.

\begin{remark} \label{freak} {\rm One may consider the problem of description of Finsler metrics whose
geodesics are circles of a fixed geodesic radius $R$ on the unit sphere.  Theorem \ref{syn} and Corollary
\ref{equall} still apply, see Remark \ref{alt}. However the situation is different on $S^2$, as far as
the continuous functions are concerned whose integrals vanish over all geodesic discs of radius $R$. For
all but countably many special values of $R$, such functions are identically zero, see \cite{Un}. This
implies an interesting ``almost everywhere" rigidity:  for a generic $R$, there is only one (up to
summation with exact 1-forms) Finsler metric whose geodesics are circles of  radius $R$; this unique
metric is an analog of the metric (\ref{const}) in $\R^2$. Of course, in the plane,  all values of the
radius $R$ are equivalent, due to similarity.}
\end{remark}

\medskip

{\bf Acknowledgments}. I am grateful to Karl Friedrich Siburg  for numerous stimulating
discussions; his visit at Penn State and a talk on the preprint \cite{P-S} were the starting points for
this work.  It is a pleasure to acknowledge fruitful discussions with J. C. Alvarez, M. Berger, D.
Khavinson, M. Levi and L. Zalcman.

\end{document}